\documentclass[12pt]{article}

\usepackage{amsfonts}
\usepackage{amssymb}
\usepackage{amsmath}

\renewcommand{\phi}{\varphi}
\newcommand{\be}{\begin{equation}}
\newcommand{\ee}{\end{equation}}

\newcommand{\nul}{{\mathbf0}}
\newcommand{\rd}{{\mathbb R}^d}

\newcommand{\zd}{{\mathbb Z}^{d}}
\renewcommand{\r}{{\mathbb R}}
\newcommand{\z} {{\mathbb Z}}

\newtheorem{theo}{Theorem}
\newtheorem{lem}[theo]{Lemma}

\newtheorem{coro}[theo]{Corollary}

\begin{document}

\title{Support of Non-separable Multivariate Scaling Function.}

\author{%
{\bf Irina Maximenko}  \thanks{Supported by RFFI, grant \# 06-01-00457.}
}
\date{}

\maketitle

\abstract{%
We make an estimation of the support of a multivariable scaling function for an arbitrary dilation matrix.
We give a method of calculating the values of the scaling function on a tight set using the knowledge
of the size of the support. 
}

\vspace{1.5cm}

\section{Introduction}
\label{ss0}

For construction of wavelet bases it is often convenient to use a scaling function, a function which
satisfies the following functional equation
\be
         \phi (x) = 2 \sum \limits_{ z \in \z}{c_{q} \phi (2x-q)},
      \ \ x \in \r.
\label{01}
\ee
The equation (\ref{01}) is called a scaling equation and 
the sequence of coefficients $\{ c_q \}$ is called a mask. The scaling function can be built from the mask.
Finite support masks are of the most interest because they generate wavelet functions with compact support. 
The knowledge of the support, particularly, allows to apply an algorithm of constructing the scaling function on a tight set.
In one-dimensional case if function $\phi$ satisfies equation
\be
         \phi (x) = 2 \sum \limits_{q=N1}^{N2} {c_{q} \phi (2x-q)},
      \ \ x \in \r,
\label{02}
\ee
then the support of scaling function is contained between N1 and N2 (you can see, for example, \cite {D1}, \cite{Ch1}).

In d-dimensional case the coefficient in scaling equation is a matrix which satisfies some natural demands.
Dilation matrix is an $ d \times d $
integer matrix such that all its eigenvalues moduluses are more than unit.
We define the norm of a matrix as the operator norm from $\rd$ to $\rd$:
$$
\|M\| = \max \sqrt {|\lambda_{MM^*}|},
$$
where $M^*$ is conjugative to $M$ matrix.
For dilation matrices we have 
\be
\lim_{n\to+\infty}\|M^{-n}\|=0  \ \ \ \ and \ \ \ \ 
\lim_{n\to+\infty} |M^{n} x|=0, \ \ \ \forall x \in \rd , \ \ \ x \ne 0 .
\label{00}
\ee
Hence the set $ \{ M^{-j} k \}_{j \in \z, \ k \in \zd} $ is tight in $\rd$.
Note that the norm of matrix $M^{-1}$ isn't obliged to be less than unit. For example matrix
$$
A=
\left(
\begin{array}{cc}
0 & 1\\
3 & 1
\end{array}
\right) 
$$
is a dilation matrix, its eigenvalues are $\approx 2.3028$ and $  -1.3028 $,
nevertheless the norm of inverse matrix approximately equals to $\approx 1.1233$, that means that $\|A^{-1}\|> 1$.
\vspace{0.3cm}

Let \ $ m:=|\det M| $ ($\det M$ is the determinant of matrix $M$), then $m$ is integer, more than unit.
The equation
\be
\label{1}
     \phi (x) = m \sum \limits_{ q \in \zd }{c_q \phi (Mx-q) } ,
      \ \ x \in \rd
\ee
is called a scaling equation and the function $ \phi $ is called $ (M,c)$ scaling function.
\vspace{0.3cm}

The cascade operator $ T_c $ is the linear operator given by
\be
     (T_c f)(x):= m \sum \limits_{ q \in \zd}    {c_q f (Mx-q)}.
\label{4}
\ee
The iteration scheme $T_c^n f=T_c(T_c^{n-1} f) , \ \ n=1,2,... $ is called a cascade algorithm. 
As a initial function $F_0$ we can take any peacewise continuous function with compact support which satisfies the following conditions
\be
\sum \limits_{q \in \zd} {F_0(x+q)}=1, \ \ \ \widehat{F}(\nul) =1 \ \ \   {\mbox{and}} \ \ \ \widehat{F}(u) \ \ \  {\mbox {is continious in the origin}}.
\label{S1a}
\ee
For example as a $F_0$ we can take the characteristic function $F_0 := \chi_{[-1/2,1/2)^d}$ or $F_0= \prod \limits_{k=1}^{d}B(x_j)$
where $B(x_j)=1-|x_j|$ if $|x_j|<1$ and $B(x_j)=0$ in another case.
\vspace{0.3cm}

Let's define 1-periodic by each variable function
$$
     m_0 (u):= \sum \limits_{q \in \zd}{c_qe^{-2 \pi i(q,u)}} , \ \ \ \ 
     \ \ u \in \rd.
$$
Obviously if $ \{ c_q \} $ has finite support then function $ m_0$  is a trigonometric polynomial.
The scaling equation in frequency domain is 
$$
\widehat{\phi}(u)=m_0({M^*}^{-1}u) \widehat {\phi}({M^*}^{-1}u). 
$$
Repeating this procedure (suppose that Fourier transform of scaling function is continuous in origin and $\widehat \phi (\nul )=1$)
we get
\be
\widehat \phi (u)=\prod \limits_{j=1}^{\infty}{m_0({M^*}^{-j}u)}.
\label{IP}
\ee
For convergence of the infinite product it is necessary that $m_0(\nul)=1$ or (the equivalent condition)
\begin{equation}
\sum \limits_{k \in \zd}{c_k}=1. 
\label{M31}
\end{equation}
If $m_0$ is a trigonometric polynomial then the infinite product converges uniformly on compact sets.

Let the space $S'$ be a set of tempered distributions defined on a set of test functions $S=S(\rd)$ (infinitely differentiable functions
which with all their derivatives decrease on infinity faster than arbitrary power function).
The sequence of distributions (or functions) $f_n$ converge in $S'$ to distribution $f$ if
$\lim \limits_{n \to \infty}(f_n, \ g) =(f, \ g)$ for all functions $g \in S$.
Distribution $f$ equals to zero in neighborhood U of the point $x_0$ if for any test function $g$ which is not zero only on 
neighborhood U we have $(f, \ g) =0$  (see for example \cite{GSH}).  
If distribution $f$ isn't equal to zero in any neighborhood U of the point $x_0$ then
$x_0$ is called essential for functional $f$. Totality of all essential points is called the support of distribution $f$.
The support of distribution $f$ corresponding to usual continuous or piecewise continuous function $f$ is the closure of the set on which
$f(x) \ne 0$.
\vspace{0.1cm}

In book \cite{NPS} it is proved that the condition (\ref{M31}) is sufficient for convergence of the cascade algorithm in $S'$. So

\begin{theo} \cite{NPS}. Let $\{ c_q \} $ be a finite mask with (\ref{M31}). Then scaling equation
(\ref{1}) has unique, up to constant  factor, decision $\phi \in S'$ with compact support. 
This decision is given by (\ref{IP}). More over, for any distribution with compact support
$f \in S'$ sequence $f_n=T^n f$ converges in $S'$ to function $c \cdot \phi$ where factor 
$c$ equals to $\widehat {f}(0)$.
\label{DSE}
\end{theo}


\vspace{0.2cm}

In section \ref{ss3} the support of a multivariate scaling function is estimated for an arbitrary dilation matrix and for masks which provide
the convergence in $S'$ of the cascade algorithm. Then obviously the estimation is true for other types of convergences.
\vspace{0.2cm}

In theory if we know the mask we can find the scaling function and the wavelet function because we can make inverse Fourier 
transform in (\ref{IP}). However we can use a simpler on practice method of construction $\phi$.
In section \ref{ss1} we show how the knowledge of size of support can be used for calculating values of a scaling function on a tight set
(in one dimensional case this method is described, for example, in  \cite{Ch1}).

\section{The Estimation of Scaling Function Support}
\label{ss3}

\begin{lem}
Let $F_1, ... , F_n, ...$ be peacewise continuous functions with compact supports on $\rd$.
Let the sequence of functions $F_1, ... , F_n, ...$ converge in $S'$. 
Suppose that $ supp \ F_1 \subset A_1 , \ supp \ F_2 \subset A_2, \ ... \ , supp \ F_n \subset A_n, \ ... $, 
where the sets $A_1, \ A_2,..., A_n, ...\in \rd$ are closed full spheres with the same center.
Let the sequence of spheres $A_n$ converge which means that there exists 
a sphere $A$ such that its radius is the limit of radiuses of spheres $A_n$. 
Then the support of limit function $F$ is in $A$. 
\label{supp}
\end{lem}

Proof.
We are going to show that any point out of $A$ doesn't belong to the support of distribution $F$.
Let's  fix an arbitrary point $a_0 \in \rd \setminus A$. Since the set $A$ is closed  there exists neighborhood $U$ of point $A$ which
is separated from the set $A$. Since $ \lim A_n =A$ there exists a number $N$ such that all spheres $A_n$, \ $n>N$ will be separated
from $U$. Let test function  $g \in S$ be not zero only in the neighborhood $U$.
Then  $(F_n , \ g)=0$ for all $n>N$ and hence $(F , \ g)=0$. 
That means that point $a_0$ doesn't belong to the support of distribution $F$.
\ \ \ \ \ \ \rule [-5pt]{5pt}{5pt}
\vspace{0.2cm}

Note. The lemma is true if instead of full spheres we consider rectangle parallelepipeds with faces parallel to coordinate hyperplanes
with the same center. The convergence of  rectangles means here the convergence of edges.
\vspace{0.3cm}

Denote $\Omega$ the set of indexes $q$ such that $c_q$ is not zero.
The mask $c_q$ has a finite support hence $\Omega$ is a finite set. Remind the definition of cascade operator
(as an initial function $F_0$ we take a piecewise continuous function with compact support which satisfies (\ref{S1a}) )
\be
F_n(x)= (T_c F_{n-1})(x):= m \sum \limits_{ q \in \Omega}   {c_{q} F_{n-1} (Mx-q)}.
\label{4C}
\ee

First we estimate the support of limit (scaling) function when $\| M^{-1} \| < 1$. 

\begin{theo}
Let $M$ be a dilation matrix and $\| M^{-1} \| < 1$. 
Suppose that mask $\{ c_q \} $ has finite support and (\ref{M31}). 
Then for the support of limit function $\phi$ the following estimation is true
$$ 
{\rm supp} \ \phi \subset \{ x \in \rd : |x| \le \frac{Q \ \| M^{-1} \|}{1 - \| M^{-1} \|} \},
$$
where $Q:=\max \limits_{q \in \Omega}|q|$.
\label{norma}
\end{theo}

Proof. Let's fix an initial function $F_0$which is a piecewise continuous function with compact support that satisfies (\ref{S1a}).
Denote its support as $\Omega_0$. The support of function $F_0$ is a compact set hence there exists a full sphere with radius R
which contains $\Omega_0$, \ $R:= \max \{ |x|, \ x \in \Omega_0 \}. $
Applying cascade operator (\ref{4C}) to function $F_0$ we get function $F_1$. Function $F_1$ is not zero if $M x-q \in \Omega_0$
for any $q \in \Omega$ or what is the same $x \in M^{-1} (\Omega_0+q)$. 
Denote  the support of function $F_1$ as $\Omega_1$. 
Then for $x \in \Omega_1$  it is true that
\be
| x | \le | M^{-1} (\Omega_0+q) | \le \| M^{-1}\| R+\| M^{-1}\| Q.
\label{Q2}
\ee
Denote the support of function $F_n$ as $\Omega_n$. Let's prove by induction that
function $F_n$ equals to zero outside of the set of $x$ where
\be
|x|  \le \| M^{-1}\|^n R+Q ( \| M^{-1}\| + \| M^{-1}\|^2+...+\| M^{-1}\|^n).
\label{Q1}
\ee
For $n=1$ the assertion is true (\ref{Q2}). Suppose that estimation  (\ref{Q1}) is true for $n$ and we prove it for $n+1$. 
Applying cascade operator (\ref{4C}) to function $F_n$ we see that function
$F_{n+1}$ is not zero if $M x-q \in \Omega_n$ for any $q$ or what is the same
$x \in M^{-1} (\Omega_n+q)$. Then for $x \in \Omega_{n+1}$ it is true that
$$
|x|  \le \| M^{-1}\| Q+ \| M^{-1}\| (\| M^{-1}\|^n R+Q ( \| M^{-1}\| + \| M^{-1}\|^2+...+\| M^{-1}\|^n))=
$$
$$
=\| M^{-1}\|^{n+1}R+Q ( \| M^{-1}\| + \| M^{-1}\|^2+...+\| M^{-1}\|^{n+1}).
$$
The induction is completed. By theorem \ref{DSE} cascade algorithm converges in $S'$.
Because $\| M^{-1}\| < 1$ the limit of the right side (\ref{Q1}) exists and is finite. 
Using lemma \ref{supp} and directing n to infinity we have
$$
|x| \le \frac{Q \ \| M^{-1} \|}{1 - \| M^{-1} \|}, \ \ \ \ t \in {\rm supp} \ \phi .
\ \ \ \ \ \ \rule [-5pt]{5pt}{5pt}
$$

Note. The support of scaling function $\phi$ doesn't depend on the size of support of the initial function $F_0$.

\begin{coro}
Let $d=1$. Then dilation matrix equals to $m, \ \ m \in \z, \ \ |m|>1$.
Suppose that mask $\{ c_q \} $ has finite support and (\ref{M31}). Then
$$ 
{\rm supp} \ \phi \subset \{ x \in \r : \ \ |x| \le \frac{Q}{|m| - 1} \}.
$$
\label{odn}
\end{coro}

\begin{coro}
Let $M$ be a diagonal dilation matrix $d \times d$.
Suppose that mask $\{ c_q \} $ has finite support and (\ref{M31}). Then
$$ 
{\rm supp} \ \phi \subset \{ x= (x_1,...,x_d) \in \rd : \ \ |x_k| \le \frac{Q}{|\lambda_k| - 1}, \ \ \ \ k=1,...,d \}.
$$
\label{diag}
\end{coro}

Let now the norm of matrix $M^{-1}$ be not necessary less than unit.
Consider first the case if matrix $M$ is an Jordan box size $s$ with eigenvalue $\lambda, \ \ |\lambda|>1$. 

\be
Mx=\left(
  \begin{array}{ccccc}
         \lambda &   0    & \ldots  &   0   &   0   \\
             1   &\lambda & \ldots  &   0   &   0   \\
             0   &   1    & \ldots  &   0   &   0   \\
         \vdots  & \vdots & \ddots  & \vdots &\vdots \\
             0   &   0    & \ldots  &   1   & \lambda \\
  \end{array}
  \right)
\left(
  \begin{array}{c}
         x_1  \\
         x_2  \\
         x_3  \\
         \vdots \\
          x_s \\
  \end{array}
  \right) =
\left(
  \begin{array}{c}
   \lambda x_1  \\
     x_1  + \lambda x_2  \\
     x_2 +  \lambda  x_3  \\
         \vdots \\
      x_{s-1} +\lambda  x_s \\
  \end{array}
  \right) .
\label{J1}
\ee

\begin{lem}
Let dilation matrix $M$ be a Jordan box size $s$ with eigenvalue $\lambda $. 
Suppose that mask $\{ c_q \} $ has finite support and (\ref{M31}). Then for each coordinate of the support of
the limit function $\phi$ we have
$$
\mbox{if}  \ \ \ \ |\lambda| \ne 2, \ \ \ \ \mbox{then} \ \ \ \  |x_k| \le \frac {Q}{|\lambda|-2} \left(1- \frac{1}{(|\lambda| - 1)^k}
\right) , \ \ \ \ k=1,...,s;
$$
$$
\mbox{if} \ \ \ \ |\lambda| = 2, \ \ \ \ \mbox{then} \ \ \ \   |x_k| \le Qk , \ \ \ \  k=1,...,s.
$$
\label{jordan}
\end{lem}
Proof. Matrix $M$ is a dilation matrix hence  $| \lambda | > 1$.
Let's fix an initial function $F_0$ a piecewise continuous function with compact support which satisfies (\ref{S1a}).
Its support $\Omega_0$ is contained in the sphere of radius R, \ $R:= \max \{ |x|, \ x \in \Omega_0 \}. $
Applying cascade operator (\ref{4C}) to function $F_0$ we get function $F_1$. Function $F_1$ is not zero if $M x-q \in \Omega_0$
for any $q \in \Omega$ or what is the same $x \in M^{-1} (\Omega_0+q)$. 
As above denote as $\Omega_1$  the support of function $F_1$. Then 
$|M x| \le R+|q| \le R +Q$ for  $ x \in \Omega_1.$
First coordinate of vector $M x$ equals to $\lambda x_1$ hence
$
|\lambda x_1 | \le |Mx| \le R +Q.
$
When we multiply matrix $M$ by vector the first coordinate doesn't depend on others.
Then by the same reasons as in theorem \ref{norma} we have
\be
|x_1^{(n)}|  \le
\frac{R}{|\lambda|^n}+Q \left( \frac{1}{|\lambda|} +
\frac{1}{|\lambda|^2}+...+\frac{1}{|\lambda|^n} \right)=:A_{n1}, \ \ 
n=1,2,3,... .
\label{suppn1}
\ee
Using lemma \ref{supp} and directing $n$ to infinity we have
\be
|x_1^{(\infty)}|  \le
Q \left( \frac{1}{|\lambda|} + \frac{1}{|\lambda|^2}+...+\frac{1}{|\lambda|^n}+...\right)=
\frac{Q}{|\lambda|-1}=:A_{\infty \ 1}.
\label{suppinf1}
\ee
Let's prove by induction that for support of function $F_1$ for each coordinate we have an estimation
\be
|x_k^{(1)}| \le \left( Q+R \right)
\left(\frac{1}{|\lambda|}+\frac{1}{|\lambda|^2}+...+
\frac{1}{|\lambda|^k}\right)=:A_{1k} , \ \ \ 
k=1,...,s.
\label{supp1k}
\ee
For $k=1$ the assertion follows from (\ref{suppn1})
$$
|x_1^{(1)}| \le \left( Q+R \right) \frac{1}{|\lambda|} .
$$
Suppose that estimation  (\ref{supp1k}) is true for $k \ (k=1,...,s-1)$ and we will prove it for $k+1$. 
$$
|x_k^{(1)}+ \lambda x_{k+1}^{(1)}| \le |Mx| \le Q+R .
$$
Using the induction postulate we have
$$
|x_{k+1}^{(1)}| \le \frac{1}{|\lambda|}
\left( Q+R + A_{1k} \right)=
$$
$$
=\frac{1}{|\lambda|}
\left( Q+R + \left( Q+R \right)
\left(\frac{1}{|\lambda|}+\frac{1}{|\lambda|^2}+...+
\frac{1}{|\lambda|^k}\right) \right)=
$$
$$
=\left( Q+R\right)
\left(\frac{1}{|\lambda|}+\frac{1}{|\lambda|^2}+...+
\frac{1}{|\lambda|^{k+1}}\right)  ,
$$
The induction is completed. So the estimation for function $F_1$  support is proved.
We get function $F_n$ from the formula (\ref{4C}) hence if $|x_k^{(n)}| \le A_{nk}$ then
$$
|x_{k-1}^{(n)}+\lambda x_k^{(n)}| \le Q+A_{n-1 \ k}
$$
Hence
\be
|x_k^{(n)}| \le \frac{1}{|\lambda|} \left( Q+A_{n-1 \ k}
+A_{n \ k-1} \right)=:A_{nk}, \ \ k=2,...,s, \ \ n=1,2,3,...
\label{suppnk}
\ee
We get a recurrent formula
$$
A_{nk}=\frac{1}{|\lambda|} \left( Q+A_{n-1 \ k}
+A_{n \ k-1} \right) .
$$
Directing $n$ to infinity we get the estimation for support of limit function $\phi$
$$
A_{\infty \ k}= \frac{1}{|\lambda|} \left( Q+A_{\infty \ k}
+A_{\infty \ k-1} \right), \ \ k=2,...,s.
$$
From this equation for  $A_{\infty \ k}$ we have
\be
A_{\infty \ k} = \frac{1}{|\lambda|-1} 
\left( Q+A_{\infty \ k-1} \right).
\label{suppinfk}
\ee
Let's prove by induction by $k$ that 
\be
A_{\infty \ k}= Q \left( 
\frac{1}{|\lambda|-1} + \frac{1}{(|\lambda|-1)^2}+...+
\frac{1}{(|\lambda|-1)^k}\right) , \ \ k=1,...,s.
\label{suppphik}
\ee
For $k=1$ the assertion is true (\ref{suppinf1}). 
Suppose now that estimation  (\ref{suppphik}) is true for $k \ (k=1,...,s-1)$, and prove it for $k+1$. 
$$
A_{\infty \ {k+1}} = \frac{1}{|\lambda|-1} 
\left( Q+A_{\infty \ k} \right)=
$$
$$
=\frac{1}{|\lambda|-1} \left( Q+Q \left( 
\frac{1}{|\lambda|-1} + \frac{1}{(|\lambda|-1)^2}+...+
\frac{1}{(|\lambda|-1)^k}\right) \right)=
$$
$$
=Q \left(  
\frac{1}{|\lambda|-1} + \frac{1}{(|\lambda|-1)^2}+...+
\frac{1}{(|\lambda|-1)^{k+1}} \right),
$$
this completes the induction. Calculating the sum of finite geometrical progression
for $|\lambda| \ne 2$ we have
\be
A_{\infty \ k}=\frac{Q}{|\lambda|-2} 
\left(1-\frac{1}{(|\lambda|-1)^k} \right) , \ \ \ k=1,...,s,
\label{otvet}
\ee
for $|\lambda| = 2$ we have
\be
A_{\infty \ k}=Qk, \ \ \ k=1,...,s.
\ \ \ \rule [-5pt]{5pt}{5pt}
\label{otvet1}
\ee

\vspace{0.1cm}

Note. By proving lemma \ref{jordan} we didn't use the fact that $q$ are integers.

\vspace{0.3cm}

An arbitrary non-degenerate matrix $M$ can be represented as $M=C^{-1}GC$, where $G$ is a Jordan matrix (consists of 
Jordan boxes), matrix C is unitary. Each of Jordan boxes (including boxes size $1 \times 1$) generates
a subspace invariant to multiplying by matrix $M$. Let's denote the rectangular parallelepiped $P \in \rd$ as follows:

\ \ if in row $p$ of matrix $G$ there is a simple eigenvalue $\lambda_p$, then
$|x_p| \le \frac{Q}{|\lambda_p| - 1}$; 

\ \ if in row $p$ of matrix $G$ there is a beginning  of a Jordan box of size $s$ which corresponds to eigenvalue
$\lambda_p \ne 2$  then $ |x_{p+k}| \le \frac {Q}{|\lambda_p|-2} \left(1- \frac{1}{(|\lambda_p| - 1)^k}
\right), \ \  k=1,...,s; $ 

\ \ if in row $p$ of matrix $G$ there is a beginning  of a Jordan box of size $s$ which corresponds to eigenvalue
$\lambda_p = 2$ then $ |x_{p+k}| \le Qk, \ \  k=1,...,s.$


\begin{theo}
Let $M$ be a dilation matrix and all its eigenvalues $\lambda_1,...,\lambda_d$ be real.
Suppose the mask $\{ c_q \} $ has finite support and (\ref{M31}). 
Then the support of limit function $\phi$ is contained in $C P$, where matrix $C$ and set $P$ are defined above.
\label{diff}
\end{theo}
Proof. The matrix $M$ has real eigenvalues hence in formula $M=CGC^{-1}$ matrices $G$ and $C$ are real.
Applying cascade operator (\ref{4C}) to function $F_n$ we get
\be
F_{n+1} (x) = m \sum \limits_{ q \in \Omega}{c_q F_n (CGC^{-1}x-q)},
      \ \ x \in \rd.
\label{s1}
\ee
Making the change of variables $x=Ct$ we get
\be
      F_{n+1} (Ct) = m \sum \limits_{ q\in \Omega}{c_q F_n (CGt-q)},
      \ \ t \in \rd .
\label{s2}
\ee
Denoting $F_n^1 (t)= F_n (Ct)$ the formula (\ref{s2}) will be:
\be
     F_{n+1}^1 (t) =
m \sum \limits_{ q\in \Omega}{c_q F_n^1 (Gt-C^{-1}q)},
      \ \ t \in \rd.
\label{s3}
\ee
Directing $n$ to infinity we have
\be
     \phi_1 (t) =m \sum \limits_{ q \in \Omega}
     {c_{q} \phi_1 (Gt-C^{-1}q)},
      \ \ t \in \rd,
\label{s5}
\ee
where $\phi_1 (t)=\phi (Ct) $. 

Matrix $G$ consists of Jordan boxes. Each of them corresponds to a subspace invariant to
multiplying by matrix $G$.
Hence using the note after lemma \ref{jordan}, we can apply lemma \ref{jordan} to function $\phi_1 $, to 
each Jordan box separately. If the corresponding eigenvalue is simple then we apply corollary \ref{odn}.
Then we get rectangular parallelepiped which we denote as $P$. 
Then we apply matrix $C$ to the bounds of support of function $\phi_1 $ and we get the set $CP$ which
will contain support of function  $\phi$.
$\ \ \ \rule [-5pt]{5pt}{5pt}$


\section{Values of Scaling Function on a Tight Set}
\label{ss1}

Let coefficients $c_q$ be not zero on the
set $\Omega$, where $\Omega$ is the finite system integer $d$-dimension vectors.
If we know the bounds of support of scaling function $\phi$, then using scaling equation (\ref{1}),
we can find values of function $\phi$ on tight set $ \{ M^{-j} k \}_{j \ge 0, \ k \in \zd} $. 
If we know the values of scaling function $\phi$ in integer points then equations
$$
     \phi (M^{-1}x) = m \sum \limits_{ q \in \Omega }{c_q \phi (x-q) } ,
$$
$$
     \phi (M^{-2}x) =m \sum \limits_{ q \in \Omega }{c_q \phi (M^{-1}x-q) } ,
$$
$$
.\ .\ .\ .\ .\ .\ .\ .\ .\ .\ .\ .\ .\ .\ .\ .\ .\ .\ .\ .\ .\ .\ .\ .
$$
uniquely define values of function $\phi$ in points $M^{-j} k, \ \ j \in \z, \ k \in \zd $.

For finding values of function $\phi$ in integer points we will use scaling equation (\ref{1}) once more.
Suppose values $\phi (k_1),...,\phi (k_N)$ are not zero. Making transformation of summation index in the right side of (\ref{1})
we get
$$
\phi (k) = m \sum \limits_{ q \in \Omega }{c_q \phi (Mk-q) } =
m \sum \limits_{ p \in \Omega_1 }{c_{Mk-p} \phi (p) }.
$$     
Denote $r=( \phi (k_1),...,\phi (k_N))^T$ and denote $B$ as the operator generated by matrix $c_{Mk-p}$. Then the matrix of
operator $B$ looks like
$$ B= m \left(
\begin{array}{cccc}
c_{Mk_1-k_1} & c_{Mk_1-k_2} & \ldots & c_{Mk_1-k_N}\\
c_{Mk_2-k_1} & c_{Mk_2-k_2} & \ldots & c_{Mk_2-k_N}\\
\vdots &  \vdots & \ddots & \vdots \\
c_{Mk_N-k_1} & c_{Mk_N-k_2} & \ldots & c_{Mk_N-k_N}
\end{array}
\right) .
$$
In the matrix form we have $r=Br.$ The values of function $\phi$ in integer points are the coordinates of eigenvector $r$
which corresponds to eigenvalue 1 of matrix $B$. The vector $r$ is normalized so that function $\phi$ satisfies partion of unit:
$$
\phi(k_1)+\phi(k_2)+...+\phi(k_N)=1.
$$
If cascade algorithm strong converges then the eigenvector corresponding to eigenvalue 1 of matrix $B$ is unique up to constant.

\end{document}